\documentclass[12pt]{amsart}
\usepackage{amsmath,amssymb,amsthm,mathrsfs}
\usepackage[vcentermath]{youngtab}
\usepackage[all]{xy}
\usepackage{amscd}
                    
\newcommand{\Sym}{\mathfrak{S}}

\newcommand{\N}{\mathbb{N}}
\newcommand{\Z}{\mathbb{Z}}
\newcommand{\Q}{\mathbb{Q}}

\newcommand{\A}{\mathcal{A}}

\newcommand{\gl}{\mathfrak{gl}}

\newcommand{\He}{\mathscr{H}}

\newcommand{\Hom}{\operatorname{Hom}}

\newcommand{\ann}{\operatorname{Ann}}

\newcommand{\Tab}{\mathrm{Tab}}

\newcommand{\dom}{\unrhd}
\newcommand{\sdom}{\rhd}
\newcommand{\notdom}{\ntrianglerighteq}
\newcommand{\lessdom}{\unlhd}
\newcommand{\slessdom}{\lhd}
\newcommand{\notlessdom}{\ntrianglelefteq}

\renewcommand{\le}{\leqslant}
\renewcommand{\ge}{\geqslant}

\newcommand{\tr}[1]{{#1^\prime}}
\newcommand\comp[1]{{{#1}^\mathrm{c}}}

\newcommand{\tilM}{\widetilde{M}}

\newcommand{\UU}{\mathbf{U}}
\newcommand{\divided}[2]{#1^{(#2)}}

\newtheorem{thm}{Theorem}%[section]
\newtheorem{lem}[thm]{Lemma}

\newtheorem{cor}[thm]{Corollary}
\newtheorem*{thm*}{Theorem}
\newtheorem*{lem*}{Lemma}
\newtheorem*{prop*}{Proposition}
\newtheorem*{cor*}{Corollary}
\theoremstyle{definition}
\newtheorem{rmk}[thm]{Remark}

\newtheorem*{rmk*}{Remark}
\newtheorem*{rmks*}{Remarks}
\newtheorem{ex}[thm]{Example}
\allowdisplaybreaks
\numberwithin{thm}{section}

\begin{document}
\title{Annihilators of permutation modules}
\author{Stephen Doty}
\address{Mathematics and Statistics,
Loyola University Chicago,
Chicago, Illinois 60626 U.S.A.}
%\curraddr{}
\email{sdoty@luc.edu}

%\thanks{}

\author{Kathryn Nyman}
\address{Mathematics Department, 
Willamette University, 900 State Street, 
Salem, Oregon 97301 U.S.A.}
\email{knyman@willamette.edu}

\date{21 April 2009}
%\dedicatory{This paper is dedicated to ...}

\begin{abstract}
Permutation modules are fundamental in the representation theory of
symmetric groups $\Sym_n$ and their corresponding Iwahori--Hecke
algebras $\He = \He(\Sym_n)$.  We find an explicit combinatorial basis
for the annihilator of a permutation module in the ``integral'' case
--- showing that it is a cell ideal in G.E.~Murphy's cell structure of
$\He$. The same result holds whenever $\He$ is semisimple, but may
fail in the non-semisimple case.
\end{abstract}

\maketitle

\section{Introduction}\noindent
Let $R$ be a commutative ring with 1 and fix an invertible element
$q\in R$. The representation theory of symmetric groups $\Sym_n$ and
the corresponding Iwahori--Hecke algebras $\He_{R,q} =
\He_{R,q}(\Sym_n)$ starts with the transitive permutation modules
$M_R^\lambda$ indexed by partitions $\lambda$. There are also twisted
versions $\tilM_R^\lambda$ and the theory may equivalently be
approached through the $\tilM_R^\lambda$ instead of the
$M_R^\lambda$. The purpose of this paper is to give explicit
combinatorial bases for the annihilators of $M_R^\lambda$ and
$\tilM_R^\lambda$ in two cases:
\begin{enumerate}\renewcommand{\theenumi}{\roman{enumi}}
\item the coefficient ring is the ring $R=\Z[v,v^{-1}]$ of ``Laurent
polynomials,'' where $v$ is an indeterminate and $q=v^2$;

\item the coefficient ring is a field $R$ such that $\He_{R,q}$ is
semisimple.
\end{enumerate}
The result is essentially the same in both cases; (i) is obtained in
Theorem \ref{thm:main-integral} by a refinement of an argument of
H\"{a}rterich \cite{Har}, and (ii) is obtained in Theorem
\ref{thm:main-ss}. It turns out that in these cases the annihilator of
$M_R^\lambda$ is a cell ideal with respect to Murphy's cellular basis
of $\He_{R,q}$. The result does not necessarily hold in case
$\He_{R,q}$ is not semisimple; see Example \ref{ex:counter}.

In case $v=1$ the algebra $\He_{R,q}$ is isomorphic to the group
algebra $R\Sym_n$ of the symmetric group $\Sym_n$; so by specializing
$v$ to 1 we obtain corresponding results on annihilators of
permutation modules for symmetric groups. Specifically, we obtain a
description of $\ann_{R\Sym_n} M_R^\lambda$ in case (i) $R = \Z$, or
(ii) $R$ is a field such that $R\Sym_n$ is semisimple.

Theorem \ref{thm:main-ss} is easily derived from a lemma of \cite{Har}
along with properties of cellular bases. To prove Theorem
\ref{thm:main-integral} we exploit Schur--Weyl duality between $\He$
and the quantized enveloping algebra $\UU(\mathfrak{gl}_n)$. This
gives an ``integral'' embedding between certain permutation modules
which seems to be new; see Lemma \ref{lem-two} for the precise
statement. This result, which may be of some interest in its own
right, is the key step in the proof of Theorem
\ref{thm:main-integral}.

\section{Symmetric groups and tableaux}\noindent
We denote by $\Sym_{S}$ the symmetric group consisting of all
bijections of a given set $S$; in particular we set $\Sym_n =
\Sym_{\{1, \,\dots, \,n\}}$.  We adopt the convention that elements of
$\Sym_n$ act on the {\em right} of their arguments, so that
compositions of permutations are read from left to right. In other
words, if $\sigma, \tau \in \Sym_n$, then $i(\sigma\tau) =
(i\sigma)\tau$, for any $i \in \{1, \dots, n\}$.

We write $\lambda \vDash n$ to indicate that $\lambda$ is a
composition of $n$, meaning that $\lambda$ is an infinite sequence
$(\lambda_1, \lambda_2, \dots)$ of nonnegative integers such that
$\sum \lambda_i = n$. The individual $\lambda_i$ are the parts of
$\lambda$, and the largest index $\ell$ such that $\lambda_j = 0$ for
all $j>\ell$ is the length, or number of parts, of $\lambda$.  Zeros
at the end of $\lambda$ are usually omitted. Any composition $\lambda$
may be sorted into a unique partition $\lambda^+$, in which the parts
are in non-strict descending order. We write $\lambda \vdash n$ to
indicate that $\lambda$ is a partition of $n$. When $\lambda \vdash
n$, we denote by $\tr{\lambda}$ the transposed partition, whose Young
diagram is obtained from the Young diagram of $\lambda$ by writing its
rows as columns.

Recall that compositions are partially ordered by dominance: given
$\lambda, \mu \vDash n$, one writes $\lambda \dom \mu$ ($\lambda$
dominates $\mu$) if $\sum_{i\le j} \lambda_i \ge \sum_{i\le j} \mu_i$
for all $j$, and one writes $\lambda \sdom \mu$ ($\lambda$ strictly
dominates $\mu$) if $\lambda \dom \mu$ and the inequality $\sum_{i\le
j} \lambda_i \ge \sum_{i\le j} \mu_i$ is strict for at least one $j$.
The notations $\mu \lessdom \lambda$, $\mu \slessdom \lambda$ are
respectively equivalent to $\lambda \dom \mu$, $\lambda \sdom \mu$.
Furthermore, we note that the dominance relation reverses when taking
transposes: 
\begin{equation}
\tr{\lambda} \lessdom \tr{\mu} \quad\text{if and only if}\quad \lambda
\dom \mu
\end{equation}
for $\lambda, \mu \vdash n$.

Let $\lambda \vDash n$. A $\lambda$-tableau $t$ is a numbering of the
boxes in the Young diagram of $\lambda$ by the numbers $1, \dots, n$
such that each number appears just once. One says that $t$ is {\em
  row-standard} if the numbers in each row are increasing read from
left to right, and {\em standard} if $t$ is row-standard and the
numbers in each column are increasing read from top to bottom.  The
group $\Sym_n$ acts naturally on the set of tableaux, on the right, as
permutations of the numbering.

\section{Murphy's bases of the Hecke algebra}\label{sec:Murphy}
\noindent
Recall the definition of $\Sym_n$ as a Coxeter group: it is the group
given by the generators $s_1, \dots, s_{n-1}$ subject to the relations

(S1) \quad $s_i^2=1$; 

(S2) \quad $s_is_js_i = s_js_is_j$ \quad if $|i-j|=1$; 
 
(S3) \quad $s_is_j = s_js_i$ \quad if $|i-j|>1$.

\noindent
We may (and will) identify the generator $s_i$ with the transposition
interchanging $i$, $i+1$ and fixing all other elements of $\{1, \dots,
n\}$. Every element of $\Sym_n$ is expressible (in many ways) as a
product of the form $w = s_{i_1} \cdots s_{i_k}$. Let $\ell(w)$, the
length of $w$, be the minimum value of $k$ in all such
expressions. Any expression of the form $w = s_{i_1} \cdots s_{i_k}$
in which $k = \ell(w)$ is \emph{reduced}.

Let $R$ be any commutative ring with 1, and fix an invertible element
$q$ of $R$.  The Iwahori--Hecke algebra, $\He_{R,q} =
\He_{R,q}(\Sym_n)$ is the associative $R$-algebra with 1 given by
generators $T_1, T_2, \ldots, T_{n-1}$ satisfying the relations

(H1) \quad $(T_i+1)(T_i-q)=0$; 

(H2) \quad $T_iT_jT_i = T_jT_iT_j$ \quad if $|i-j|=1$; 
 
(H3) \quad $T_iT_j = T_jT_i$ \quad if $|i-j|>1$.

\noindent
In order to simplify notation, we shall henceforth write $\He_R$
instead of $\He_{R,q}$, letting $q$ be understood.

\begin{rmk}  
  Some authors use a different, but equivalent, quadratic relation in
  place of (H1): $(T_i - v)(T_i +v^{-1})=0$.  Setting $v=q^{\frac12}$
  and extending the scalars to include $q^{\frac12}$, one can see that
  this version is isomorphic to the one defined above. The results of
  this paper hold for the alternative version of $\He_{R}$, although
  many of the specific formulas are somewhat different.
\end{rmk}

If $w = s_{i_1} \cdots s_{i_k}$ is a reduced expression for $w\in
\Sym_n$, one defines
\[
T_w = T_{i_1} \cdots T_{i_k}.
\]
In particular, $T_{\mathit{id}}=1$ and $T_{s_i} = T_i$. The element
$T_w$ is well defined independently of the choice of reduced
expression for $w$ and satisfies
\begin{equation}
  T_w T_{s_i} = 
  \begin{cases}
   T_{ws_i} & \text{ if } \ell(ws_i) > \ell(w);\\
   qT_{ws_i} + (q-1)T_w &  \text{ if } \ell(ws_i) < \ell(w).
  \end{cases}
\end{equation}
The elements $T_w$, for $w \in \Sym_n$, form an $R$-basis of $\He_R$; in
particular $\He_R$ is free as an $R$-module. The generator $T_i=T_{s_i}$
is invertible in $\He_R$, with
\begin{equation}
  T_i^{-1} = q^{-1} (T_i - q+1) 
\end{equation}
for all $i = 1, \dots, n-1$. It follows that any $T_w$ is invertible,
with 
\begin{equation}
  T_w^{-1} = T_{i_k}^{-1} \cdots T_{i_1}^{-1} 
\end{equation}
where $w = s_{i_1}\cdots s_{i_k}$ is a reduced expression for $w$.

We refer the reader to \cite[Lemma 2.3]{Murphy2} for a proof of the
following well known result.

\begin{lem}
  Suppose $R$ is an integral domain.  Let $*$, $\dag$, $\sharp$ be the
  maps $\He_R \to \He_R$ defined on basis elements by the rules
  
  \qquad $*: T_w \to T_{w^{-1}}$; 

  \qquad $\dag: T_w \to (-q)^{\ell(w)} T_w^{-1}$;

  \qquad $\sharp: T_w \to (-q)^{\ell(w)} (T_{w^{-1}})^{-1}$

  \noindent
  for $w \in \Sym_n$, extended to $\He_R$ by linearity. Then $*$ and
  $\dag$ are commuting anti-involutions of $R$-algebras. The map
  $\sharp$ is the composite of $*$ and $\dag$; thus $\sharp$ is an
  involution of $R$-algebras.
\end{lem}

We denote the image of any $h \in \He_R$ under the maps $*$, $\dag$, and
$\sharp$ by $h^*$, $h^\dag$, and $h^\sharp$ respectively.

We now recall some results of Murphy \cite{Murphy, Murphy2}, which
provide two cellular bases of $\He_R$.  Let $\lambda \vDash n$, and
let $t^\lambda$ be the tableau of shape $\lambda$ in which the numbers
$1, \dots, n$ have been inserted in the boxes in order from left to
right along the rows.  Let $\Sym_\lambda$ be the
row stabilizer of $t^\lambda$. Then
\begin{equation}
\Sym_\lambda = \Sym_{ \{1, \dots, \lambda_1\} } \times \Sym_{
\{\lambda_1+1, \dots, \lambda_1+\lambda_2\} } \times \cdots
\end{equation}
is a Young subgroup of $\Sym_n$. Define elements
\begin{equation}
  \textstyle x_\lambda = \sum_{w \in \Sym_\lambda} T_w; \quad
  y_\lambda = \sum_{w \in \Sym_\lambda} (-q)^{-\ell(w)} T_w .
\end{equation}
Given a tableau $t$ of shape $\lambda$, for $\lambda \vDash n$, let
$d(t)$ be the unique element of $\Sym_n$ such that $t = t^\lambda
d(t)$.  Given any pair $s,t$ of row-standard $\lambda$-tableaux,
following Murphy we set
\begin{equation}
 x_{st} = T_{d(s)}^* x_\lambda\, T_{d(t)}, \quad y_{st} = T_{d(s)}^*
  y_\lambda\, T_{d(t)}.
\end{equation}
Since $x_\lambda$ and $y_\lambda$ are invariant under $*$, it follows
that
\begin{equation} \label{eq:star}
x_{st}^{*} = x_{ts}^{\ }; \quad  y_{st}^{*} = y_{ts}
\end{equation}
for any pair $s,t$ of row-standard $\lambda$-tableaux.  For any
$\lambda \vdash n$ let $\Tab(\lambda)$ be the set of standard
$\lambda$-tableaux, and set
\begin{equation*}
\He_R[\sdom \lambda] = \sum_{a,b \in \Tab(\mu),\, \mu \sdom \lambda}
Rx_{ab}; \quad \He_R[\dom \lambda] = \sum_{a,b \in \Tab(\mu),\, \mu \dom
  \lambda} Rx_{ab}.
\end{equation*}
It follows from the first equality in \eqref{eq:star} and part (b) of
the next result that $\He_R[\sdom \lambda]$ and $\He_R[\dom \lambda]$ are
both two-sided ideals of $\He_R$.

\begin{thm}\cite[Theorem 4.17 and Theorem 5.1]{Murphy2} \label{thm:Murphy} 
Assume that $R$ is an integral domain.
\item(a) The set $\{x_{st}: s,t \in \Tab(\lambda), \lambda \vdash n
  \}$ is an $R$-basis of $\He_R$.

\item(b) If $\lambda \vdash n$, for any $h \in \He_R$, $s,t \in
  \Tab(\lambda)$ we have
\[
  x_{st} h = \sum_{u \in \Tab(\lambda)} r_h(t,u) x_{su}
  \pmod{\He_R[\sdom \lambda]} 
\]
where $r_h(t,u) \in R$ is independent of $s$. 
\end{thm}

Note that by applying $*$ the equality in the theorem may be written
in the equivalent form 
\begin{equation*}
   h^* x_{ts}  = \sum_{u \in \Tab(\lambda)} r_h(t,u) x_{us}
  \pmod{\He_R[\sdom \lambda]}. 
\end{equation*}
This is used, for instance, in proving that $\He_R[\sdom \lambda]$ and
$\He_R[\dom \lambda]$ are two-sided ideals of $\He_R$.

\begin{rmk}
\item(a) The basis in part (a) of the theorem is a cellular basis in
  the sense of \cite{GL}. This follows from the first equality in
  \eqref{eq:star} and part (b) of the theorem.
\item(b) For the sake of completeness, we mention that the theorem
  remains true if one replaces $x_{st}$ by $y_{st}$ throughout. Thus
  $$\{y_{st}: s,t \in \Tab(\lambda), \lambda \vdash n \}$$ is another
  cellular basis of $\He_R$. We will not need this fact in the paper.
\item(c) By applying the involution $\sharp$ to any cellular basis of
  $\He_R$, we obtain another cellular basis. Thus $$\{x_{st}^\sharp: s,t
  \in \Tab(\lambda), \lambda \vdash n \}, \quad \{y_{st}^\sharp: s,t
  \in \Tab(\lambda), \lambda \vdash n \}$$ are both cellular bases of
  $\He_R$.
 %\item(d) CHECK RELATIONSHIP $x^{\sharp}_{st}$ and $y_{st}$.
\end{rmk}

For $\lambda \vdash n$, following \cite{GL} we let $C_R(\lambda)$ be
the free $R$-module with basis $\{ c_t: t \in \Tab(\lambda) \}$. We
define an action of $\He_R$ on $C_R(\lambda)$ by the rule
\begin{equation}
  c_t h = \sum_{u \in \Tab(\lambda)} r_h(t,u) c_u
\end{equation}
for any $t \in \Tab(\lambda)$. The element $r_h(t,u) \in R$ is
determined as in part (b) of Theorem \ref{thm:Murphy}. The
$\He_R$-modules $C_R(\lambda)$ are known as \emph{cell modules} in the
terminology of \cite{GL}. Note that for any given fixed $s \in
\Tab(\lambda)$, the right $\He_R$-module $C_R(\lambda)$ is isomorphic
with the submodule of $\He_R[\dom \lambda]/\He_R[\sdom \lambda]$
spanned by the set of all right cosets of the form
\[
  x_{st} + \He_R[\sdom \lambda] \qquad (t \in \Tab(\lambda)).
\]
Murphy has identified the cell modules for $\He_R$ with respect to the
cellular basis of Theorem \ref{thm:Murphy}. In \cite[Theorem
  5.3]{Murphy2} he proves:
\begin{equation}
  C_R(\lambda) \simeq S_{\lambda,R} \quad\text{for any } \lambda \vdash n.
\end{equation}
Here $S_{\lambda,R}$ is the linear dual $\Hom_R(S_R^\lambda, R)$ of
the Specht module $S_R^\lambda$.

The module $S_R^\lambda$, which is a $q$-analogue of the corresponding
Specht module for the symmetric group $\Sym_n$, was introduced in
\cite[\S4]{DJ:Hecke}.  Although we need only generic properties of
Specht modules associated to their interpretation as cell modules,
some readers may prefer an explicit construction. One approach is to
define $S_R^\lambda$ as the following right ideal of $\He_R$:
\begin{equation}
S_R^\lambda:= x_\lambda T_{w_\lambda} y_{\tr{\lambda}}
\He_R,
\end{equation}
where the element $w_\lambda$ is the unique element of $\Sym_n$ such
that $t^\lambda w_\lambda = t_\lambda$. Here $t_\lambda$ is the
tableau of shape $\lambda$ in which the numbers $1, \dots, n$ have
been inserted in the boxes in order from top to bottom in the columns.

\section{Permutation modules}\label{sec:perm}\noindent
In \cite{DJ:Hecke}, \cite{DJ:qsa}, Dipper and James studied the right
ideals $M_R^\lambda := x_\lambda \He_R$, for $\lambda \vDash n$,
noting that they are $q$-analogues of the classical permutation
modules (see \cite{James}) for symmetric groups. They also studied the
right ideals $\tilM_R^\lambda := y_\lambda \He_R$; these are
$q$-analogues of the signed permutation modules for symmetric
groups. We wish to study the annihilators of these $\He_R$-modules.

By \cite[(2.1)]{DJ:qsa}, there exist elements $r_\lambda$, $r'_\lambda
\in R$ such that $r_\lambda r'_\lambda = 1$ and
\begin{equation} \label{eq:xysharp}
    x_\lambda^\sharp = r_\lambda y_\lambda; \quad y_\lambda^\sharp =
    r'_\lambda x_\lambda.
\end{equation}
Given any right $\He_R$-module $N$, we obtain a new right
$\He_R$-module $N^\sharp$ in the usual way, by letting $N^\sharp = N$
and twisting the original $\He_R$-action by the automorphism
$\sharp$. On the other hand, if $N$ is a right ideal of $\He_R$ then
$N^\sharp = \{ n^\sharp: n \in N\}$ is another right ideal of
$\He_R$. Thus, in the case when $N$ is a right ideal, the notation
$N^\sharp$ has two meanings. However, the reader can easily check that
the two possible interpretations lead to isomorphic right
$\He_R$-modules.  From \eqref{eq:xysharp}, it follows immediately that
as right $\He_R$-modules we have an isomorphism
\begin{equation}
  (M_R^\lambda)^\sharp \simeq \tilM_R^\lambda
\end{equation}
for any $\lambda \vDash n$. Thus, any description of $\ann_{\He_R}
M_R^\lambda$ will give immediately a description of $\ann_{\He_R}
\tilM_R^\lambda$, simply by applying the involution $\sharp$. So we
focus on obtaining a description of the former.

Let $\He_R(\lambda)$ be the subalgebra of $\He_R$ generated by all the
$T_i$ except $T_{\lambda_1}$, $T_{\lambda_1+\lambda_2}$, \dots,
$T_{\lambda_1+\cdots+\lambda_k}$ where $\lambda = (\lambda_1,
\lambda_2, \dots, \lambda_k)$ has length $k$.  By \cite[Corollary
1.14]{Mathas}, there are two one-dimensional representations
$1_{\He_R}$ and $\varepsilon_{\He_R}$ of $\He_R$ defined on basis
elements by
\begin{equation}
  1_{\He_R}(T_w) = q^{\ell(w)}; \quad \varepsilon_{\He_R}(T_w) =
  (-1)^{\ell(w)} 
\end{equation}
for any $w \in \Sym_n$. These are known as the \emph{trivial} and
\emph{sign} representations, respectively; by abuse of notation we
denote the corresponding right $\He_R$-modules by the same symbols.
As right $\He_R$-modules, we have isomorphisms
\begin{equation}
M_R^\lambda \simeq 1_{\He_R} \otimes_{\He_R(\lambda)} \He_R, \quad
\tilM_R^{\lambda} \simeq \varepsilon_{\He_R} \otimes_{\He_R(\lambda)}
\He_R
\end{equation}
thus justifying the terminology ``permutation'' and ``signed
permutation'' for these modules.

Let $\lambda \vDash n$.  For typographical reasons, we denote by
$x_{\lambda t}$ the element $x_{st}$ with $s = t^\lambda$. Thus
$x_{\lambda t} = x_\lambda T_{d(t)}$, for any row-standard tableau $t$
of shape $\lambda$. We will need the following result of Dipper and
James, which gives a basis of $M_R^\lambda$ and determines the action of
$\He_R$ on basis elements.

\begin{lem}\cite[Lemma~3.2]{DJ:Hecke} \label{lem:Mbasis} Assume
  that $R$ is an integral domain, and let $\lambda \vDash n$.

\item(a) $\{ x_{\lambda t} : t \text{ row-standard of shape $\lambda$}
  \}$ is an $R$-basis of $M_R^\lambda = x_\lambda \He_R$. 

\item(b) Suppose $t$ is row-standard of shape $\lambda$, and set
  $u=ts_i$. Let $\mathrm{row}_t(j)$ be the row index of $j$ in
  $t$. Then
\[
x_{\lambda t} T_i = 
 \begin{cases}
   q x_{\lambda t} & \text{if } \mathrm{row}_t(i) = \mathrm{row}_t(i+1)\\
   x_{\lambda u} & \text{if } \mathrm{row}_t(i) < \mathrm{row}_t(i+1)\\
   q x_{\lambda u} + (q-1) x_{\lambda t} & \text{if } \mathrm{row}_t(i) 
     > \mathrm{row}_t(i+1).
 \end{cases}
\]
\end{lem}

Note that for $\lambda \vDash n$, $M_R^\lambda \simeq M_R^{\lambda^+}$,
where $\lambda^+$ is the unique partition in the $\Sym_n$-orbit of
$\lambda$.  Thus, when considering the annihilator of $M_R^\lambda$, it
is enough to restrict our attention to the case in which $\lambda
\vdash n$.

\section{The annihilator in the semisimple case}\label{sec:main}
\noindent
We will require a preparatory lemma of H\"arterich.  To formulate it,
we need to extend the dominance order on compositions to the set of
row-standard tableaux, as follows. Let $t$ be a row-standard
$\lambda$-tableau, where $\lambda \vDash n$. For any $j \le n$ denote
by $t_{\downarrow j}$ the row-standard tableau that results from
throwing away all boxes of $t$ containing a number bigger than
$j$. Let $[t_{\downarrow j}]$ be the corresponding composition of $j$
(the composition defining the shape of $t_{\downarrow j}$). Given
row-standard tableaux $s$ and $t$ with the same number $n$ of boxes,
define

 $s \dom t$ if for each $j \le r$, $[s_{\downarrow j}] \dom
 [t_{\downarrow j}]$;

 $s \sdom t$ if for each $j \le r$, $[s_{\downarrow j}] \sdom
 [t_{\downarrow j}]$.

\noindent
Note that if $s$ and $t$ are standard tableaux, respectively of shape
$\lambda$ and $\mu$, where $\lambda$ and $\mu$ are partitions of $n$,
then $s \dom t$ if and only if $\tr{t} \dom \tr{s}$. Here $\tr{t}$
denotes the transposed tableau of $t$, obtained from $t$ by writing
its rows as columns.  The dominance order on tableaux extends
naturally to pairs of tableaux, by defining:
\begin{equation}
    (s,t) \dom (u,w) \text{ iff } s\dom u \text{ and } t\dom w.
\end{equation}

For a set $P$ of partitions of $n$, set $\He_R[P] = \sum_{s,t} Rx_{st}$,
where the sum is taken over the set of pairs $s,t \in \Tab(\mu)$, for
$\mu \in P$. If $P$ is closed with respect to $\dom$ (i.e., $\nu
\vdash n$, $\mu \in P$, and $\nu \dom \mu$ implies $\nu \in P$), then
$\He_R[P]$ is a two sided ideal of $\He_R$.  For $\lambda \vdash n$, the
set $P = \{ \mu \vdash n: \mu \notlessdom \tr{\lambda} \}$ is closed
with respect to $\dom$, so $\He_R[\notlessdom \tr{\lambda}]$ is a
two-sided ideal of $\He_R$. Thus 
\[
 \He_R[\notlessdom \tr{\lambda}]^\sharp = \{h^\sharp | h \in
 \He_R[\notlessdom \tr{\lambda}]\}
\]
is also a two-sided ideal of $\He_R$. We remark that these ideals are
free as $R$-modules and have $R$-rank $\textstyle \sum_{\tr{\mu}
  \notlessdom \tr{\lambda}} (\dim_R S_R^{\mu})^2$. This follows from
the cellular axioms and the various identifications made above.  

The following result is a variant of \cite[Lemma 3]{Har}.

\begin{lem}[H\"arterich]\label{lem-one}
 Let $R$ be an integral domain.  For any $\lambda \vdash n$ we have
 $\He_R[\notlessdom \tr{\lambda}]^\sharp \subseteq \ann_{\He_R}
 M_R^\lambda$.
\end{lem}

\begin{proof}
  Let $a,b \in \Tab(\mu)$ where $\mu \notlessdom \lambda'$, and let
  $c$ be any row-standard $\lambda$-tableau.  We have $c \notlessdom
  a'$, since otherwise $\lambda =[c_{\downarrow n}] \lessdom
  [a'_{\downarrow n}]= \mu '$.  By Lemma 4.12 of \cite{Murphy2},
  $x^{\dag}_{ab}x_{c\lambda} = 0$.  Applying $*$, we have
  $x^*_{c\lambda}(x^{\dag}_{ab})^* = x_{\lambda
    c}x^\sharp_{ab}= 0$.  Now, since $\left \{x_{\lambda c} : c \text{ is
    a row-standard $\lambda$-tableau} \right\}$ is a basis of
  $M_R^{\lambda}$, we have $x^\sharp_{ab} \in \ann_{\He_R} M_R^{\lambda}$ for
  $a,b$ of shape $\mu \notlessdom \lambda'$.
\end{proof}

\begin{thm} \label{thm:main-ss}
Let $R$ be a field, and suppose that $\He_R$ is semisimple. Then, for
any $\lambda \vdash n$, equality holds in the inclusion in the
preceding lemma; i.e.,\ we have $\He_R[\notlessdom
  \tr{\lambda}]^\sharp = \ann_{\He_R} M_R^\lambda$.
\end{thm}

\begin{proof}
The semisimplicity of $\He_R$ implies that $M_R^\lambda$ is completely
reducible as an $\He_R$-module, with irreducible factors of the form
$S_R^\mu$ for various $\mu \dom \lambda$. By the $q$-analogue of Young's
rule
\cite[Theorem 7.2]{Murphy2}, each $S_R^\mu$ for any $\mu$ satisfying
$\mu \dom \lambda$ occurs at least once in $M_R^\lambda$. Thus, it
follows by the theory of semisimple algebras that
\[
\dim_R (\He_R/ \ann_{\He_R} M_R^\lambda) \ge \textstyle \sum_{\mu \dom
  \lambda} (\dim_R S_R^\mu)^2,
\]
or, equivalently,
\[
\dim_R \ann_{\He_R} M_R^\lambda \le \textstyle \sum_{\mu \notdom \lambda}
(\dim_R S_R^\mu)^2.
\]
But the remark at the end of the paragraph preceding Lemma
\ref{lem-one} gives the equality
\[
\textstyle \sum_{\tr{\mu} \notlessdom \tr{\lambda}} (\dim_R
S_R^{\mu})^2 = \dim_R \He_R[\notlessdom \tr{\lambda}]^\sharp
\]
and applying the equivalence between the two conditions $\tr{\mu}
\notlessdom \tr{\lambda}$ and $\mu \notdom \lambda$, we conclude that
\[
\dim_R \ann_{\He_R} M_R^\lambda \le \dim_R \He_R[\notlessdom
\tr{\lambda}]^\sharp .
\]
Finally, Lemma \ref{lem-one} gives the opposite inequality, thus
proving the result.
\end{proof}

\begin{rmk}
  Let $R$ be a field. It is well known \cite[Theorem 4.3]{DJ4} that
  $\He_R$ is semisimple unless one of the following holds:
  \begin{enumerate}
  \item $q \ne 1$ and $q$ is a primitive $e$th root of 1, where $e \le
    n$;
  \item $q=1$ and the characteristic of $R$ is $\le n$. 
  \end{enumerate}
\end{rmk}

In case $\He_R$ is not semisimple, the inclusion of Lemma \ref{lem-one}
may be strict (see Example~\ref{ex:counter}).

\section{Tensor space}\label{sec:tensor}\noindent
Henceforth we will work over the ring $\A = \Z[v, v^{-1}]$ where $v$
is an indeterminate, with quotient field $\Q(v)$. We set $q=v^2$ and
consider the Hecke algebra $\He_\A$ over the ring $\A$, defined by the
generators and relations as in \S\ref{sec:Murphy}.  For any
commutative ring $R$ with a chosen invertible element $v$ we set
$q=v^2$ in $R$. Then
\begin{equation}\label{eq:basechgiso}
  R \otimes_\A (\He_\A) \simeq \He_R,
\end{equation}
where $R$ is regarded as an $\A$-algebra by means of the ring
homomorphism $\A \to R$ sending $v\in \A$ to $v\in R$. The isomorphism
in \eqref{eq:basechgiso} is determined by sending $1\otimes T_i \to
T_i$ for all $i = 1, \dots, n-1$. Moreover, we have
\begin{equation}\label{eq:Mbasechgiso}
  R \otimes_\A (M_\A^\lambda) \simeq M_R^\lambda
\end{equation}
for any $\lambda \vDash n$, and any $R$. Our goal is to determine the
annihilator of the integral permutation module $M_\A^\lambda$.

We will need the Drinfeld--Jimbo quantized enveloping algebra
corresponding to the Lie algebra $\gl_m$ for $m \ge 2$. Let $Y$ be the
free abelian group with basis $H_1, \ldots, H_m$.  Let $\varepsilon_1,
\ldots, \varepsilon_m \in X:=Y^*$ be the corresponding dual basis;
$\varepsilon_i$ is given by $\varepsilon_i(H_j):=\delta_{i,j}$ for
$j=1, \ldots, m$. For $i=1,\ldots,m-1$ let $\alpha_i\in X$ be given by
$\alpha_i := \varepsilon_i - \varepsilon_{i+1}$. Define a partial
order $\le$ on $X$ by $\lambda \le \mu$ if and only if $\mu-\lambda
\in \sum_i \N\alpha_i$.

Define an associative algebra $\UU = \UU(\gl_n)$ (with 1) over $\Q(v)$
by the generators

\qquad $E_i, F_i \quad(i=1, \ldots, m-1), \qquad v^{h} \quad(h\in Y)$ 

\noindent
subject to the defining relations

    (U1)\quad $v^0=1$, \quad $v^hv^{h'}=v^{h+h'}$ 

    (U2)\quad $v^hE_iv^{-h}=v^{\alpha_i(h)}E_i$, \quad
        $v^hF_iv^{-h}=v^{-\alpha_i(h)}F_i$ 

    (U3)\quad $E_iF_j-F_jE_i=\delta_{ij}\frac{K_i-K_i^{-1}}{v-v^{-1}}$ \quad
    where $K_i := v^{H_{i}-H_{i+1}}$ 

    (U4)\quad $E_i^2E_j-(v+v^{-1})E_iE_jE_i+E_jE_i^2=0$ \quad if $|i-j|=1$

    (U5)\quad $F_i^2F_j-(v+v^{-1})F_iF_jF_i+F_jF_i^2=0$ \quad if $|i-j|=1$

    (U6)\quad $E_iE_j=E_jE_i,\quad F_iF_j=F_jF_i$ \quad if $|i-j|>1$ 

\noindent
for $1\le i,j \le m-1$ and $h, h' \in Y$.  The subalgebra of $\UU$
generated by all the $E_i, F_i, K_i$ (for $i = 1, \dots, m-1$) is
denoted by $\UU(\mathfrak{sl}_m)$; this is the quantized enveloping
algebra corresponding to the Lie algebra $\mathfrak{sl}_m$. 

In the remainder of this section all tensor products will be over
$\Q(v)$ unless specified otherwise.  There exist unique algebra maps
$\Delta:\UU \to \UU\otimes \UU$ (where $\UU\otimes \UU$ is regarded as
an algebra in the usual way) and $\epsilon:\UU \to \Q(v)$ such that
\begin{gather*}
\Delta(v^h)=v^h\otimes v^h,\quad \epsilon(v^h)=1\\
\Delta(E_i)=E_i\otimes K_i + 1 \otimes E_i,\quad
\epsilon(E_i)=0\\
\Delta(F_i)=F_i\otimes 1 + K_i^{-1} \otimes F_i, \quad
\epsilon(F_i)=0
\end{gather*}
for any $i = 1, \dots, m-1$, and $h \in Y$.  The map $\Delta$ defines a
comultiplication and the map $\epsilon$ a counit which together define
a coalgebra structure on $\UU$. (Actually, there is a well defined
antipode $S: \UU \to \UU^{\mathit{opp}}$ which together with $\Delta$,
$\epsilon$ give $\UU$ a Hopf algebra structure, but we shall not need
it.) The map $\Delta$ induces a map $\divided{\Delta}{n}: \UU \to
\UU^{\otimes n}$ defined as the composite 
\begin{equation*}
  \begin{CD}
    \UU @>\Delta>> \UU \otimes \UU @>\Delta\otimes 1>> \UU \otimes \UU
    \otimes \UU @>\Delta\otimes 1 \otimes 1>> \cdots @>\Delta\otimes 1
    \otimes \cdots \otimes 1>> \UU^{\otimes n} 
  \end{CD}
\end{equation*}
for $n \ge 2$.  Note that $(\Delta \otimes 1)\Delta = (1 \otimes
\Delta)\Delta$, etc.  We have $\Delta = \Delta^{(2)}$ and
\begin{gather*}
\Delta^{(n)}(v^h) = v^h \otimes v^h \otimes \cdots \otimes
v^h;\\ 
\Delta^{(n)}(E_i) = E_i\otimes K_i \otimes \cdots \otimes K_i\ +\ 
1 \otimes E_i \otimes K_i \otimes \cdots \otimes K_i\ +\\ 
\cdots\ +\ 1 \otimes \cdots \otimes 1 \otimes E_i;\\
\Delta^{(n)}(F_i) = F_i\otimes 1 \otimes \cdots \otimes 1\ +\ 
K_i^{-1} \otimes F_i \otimes 1 \otimes \cdots \otimes 1\ +\\ 
\cdots\ +\ K_i^{-1} \otimes \cdots \otimes K_i^{-1} \otimes F_i.
\end{gather*}
The map $\Delta^{(n)}$ is used to put a $\UU$-module structure on the
$n$th tensor power of a given $\UU$-module.

Set $V = \Q(v)^m$ with canonical basis $\{e_1, \dots, e_m\}$. Define a
left action of $\UU$ on $V$ by
\begin{equation}
  E_i e_j = \delta_{j,i+1} e_i; \quad F_i e_j = \delta_{j,i} e_{i+1};
  \quad v^h e_j = v^{\varepsilon_j(h)} e_j
\end{equation}
for $j=1, \dots, m$, $i=1, \dots, m-1$, $h \in Y$. This action makes
$V$ into a $\UU$-module, and by applying $\divided{\Delta}{n}$ one
makes the tensor power $V^{\otimes n}$ into a $\UU$-module as well. 

Let $I(m,n)$ be the set of all finite sequences $\mathbf{i} = (i_1,
\dots, i_n)$ where each $i_j \in \{1, \dots, m\}$. Set
\begin{equation}\label{eq:tensor}
  e_\mathbf{i}:= e_{i_1}\otimes e_{i_2} \otimes \cdots \otimes e_{i_n};
\end{equation}
in terms of this notation $\{ e_\mathbf{i}: \mathbf{i} \in I(m,n) \}$
is a basis for $V^{\otimes n}$. Let $\UU^0$ be the subalgebra of $\UU$
generated by the $v^h$ for $h \in Y$. Since
\begin{equation}
  v^h e_\mathbf{i} = v^{\lambda(h)} e_\mathbf{i}
\end{equation}
for all $h \in Y$, the vectors $e_{\mathbf{i}}$ are weight vectors (of
weight $\lambda$) for the action of $\UU^0$. Here $\lambda =
\lambda_1\varepsilon_1 + \cdots + \lambda_m\varepsilon_m \in X$ where
each $\lambda_k = $ the number of $j$ such that $i_j = k$. It follows
that
\begin{equation} \label{eq:wt-space}
  V^{\otimes n} = \textstyle \bigoplus_{\lambda \in X} (V^{\otimes n})_\lambda
\end{equation}
where $(V^{\otimes n})_\lambda$ is the $\Q(v)$-span of all
$e_\mathbf{i}$ of weight $\lambda$. 

The algebra $\UU$ has $q$-analogues of the standard root vectors in
$\mathfrak{gl}_m$, defined as follows. For a positive root
$\alpha=\varepsilon_i - \varepsilon_j$ for $1 \le i<j \le m$, if $j-i
= 1$ then set $E_{i,j} = E_i$ and $E_{j,i} = F_i$. If $j-i > 1$, we
assume by induction that $E_{i+1,j}$ and $E_{j,i+1}$ have already been
defined, and set (following Jimbo \cite[Proposition 1]{Jimbo}, Xi
\cite[Section 5.6]{Xi})
\begin{equation}
E_{i,j} = v^{-1} E_i E_{i+1,j} - E_{i+1,j} E_i; \quad
E_{j,i} = v E_{j,i+1} F_i - F_i E_{j,i+1}.
\label{eq:5.7}
\end{equation}
Set $K_{i,j} = v^{H_i-H_j}$ for each $1 \le i<j \le m$. We note the
following results.

\begin{lem} \label{lem:sl2triple}
  For $1 \le i < j \le m$ the triple $\{ E_{i,j}, K_{i,j}, E_{j,i} \}$
  generates a subalgebra of $\UU$ isomorphic with
  $\UU(\mathfrak{sl}_2)$. 
\end{lem}

\begin{proof}
 Set $E=E_{i,j}$, $K=K_{i,j}$, and $F=E_{j,i}$ for $i<j$.  We need
 only show that $E$, $K$, and $F$ satisfy the defining relations for
 $\UU(\mathfrak{sl}_2)$:
 
 (a) $KEK^{-1} = v^2 E$; \qquad $KFK^{-1} = v^{-2}F$
 
 (b) $EF-FE = \frac{K-K^{-1}}{v-v^{-1}}$
 
 \noindent In the case $j-i=1$, (a) follows from (U1) and (U2) while
 (b) is immediate from (U3).  Assume $j-i >1$, and assume that (a) and
 (b) hold for $E'=E_{i+1,j}$, $K'=K_{i-1,j}$, and $F'=E_{j, i+1}$ in
 place of $E$, $K$, and $F$ respectively.  Then, using \ref{eq:5.7}
 and the equality $K=K_iK'$, we have
 \begin{eqnarray*}
 KEK^{-1} & = & (K_iK')(v^{-1}E_iE' - E'E_i)(K_i K')^{-1}\\
 & = & v^{-1}(K_iE_iK_i^{-1})(K'E'K'^{-1})- (K'E'K'^{-1})(K_iE_iK_i^{-1}) \\
 &=& v^2 (v^{-1}E_iE' - E'E_i)\\
 &=& v^2 E
 \end{eqnarray*}
 which proves the first relation in (a).  The second relation in (a)
 and relation (b) are proved by similar calculations.
\end{proof}

The following result was observed in \cite[Lemma 4.2]{Tony} although
in a slightly different context. We adapt the argument from
\cite{Tony} to our situation.

\begin{lem}\label{lem:Tony}
Let $\mu \in X$, and set $\mu = \sum_k \mu_k \varepsilon_k$.  Assume
$\mu_k \geq 0$ for $k =1, \ldots, m$.  Suppose that $i<j$ and $\mu_i -
\mu_j >0$.  Then the action of $E_{j,i}$ gives a $\Q(v)$-linear
injection $(V^{\otimes n})_\mu \to (V^{\otimes
  n})_{\mu-\varepsilon_i+\varepsilon_j}$.
\end{lem}

\begin{proof}
It is clear that the action of $E_{j,i}$ gives a linear map from
$(V^{\otimes n})_\mu$ into $(V^{\otimes
  n})_{\mu-\varepsilon_i+\varepsilon_j}$, so we need only to prove
that the map is injective.  Set
\[ E=E_{i,j}, \qquad   K=K_{i,j}, \qquad F=E_{j,i}. \]
By induction on $r$, one proves easily that
\begin{equation*}
E^rF -FE^r = \mathop{\sum_{s+s' = r-1}}_{ s,s' \geq 0} E^s
\frac{K-K^{-1}}{v-v^{-1}} E^{s'}
\end{equation*}
for any positive integer $r$.  Fix some $0 \neq u \in (V^{\otimes
  n})_\mu$ and choose $r$ so that $E^r u=0$ but $E^s u \neq 0$ for any
$s < r$.  This is possible since $E$ acts locally nilpotently on any
finite-dimensional $\UU$-module.  Then
\begin{eqnarray*}
E^rF u & = & (E^rF - FE^r)u \\ & = & \mathop{\sum_{s+s' = r-1}}_{ s,s'
  \geq 0} E^s \frac{K-K^{-1}}{v-v^{-1}} E^{s'} u.
\end{eqnarray*}
But a calculation in $\UU(\mathfrak{sl}_2)$ shows that 
\[ \frac{K-K^{-1}}{v-v^{-1}} E^{s'} u = v^{\mu_i - \mu_j + 2s'} E^{s'} u, \]
and it follows that $E^r F u = \left( \sum_{s' = 0}^{r-1} v^{\mu_i -
  \mu_j + 2s'} \right) E^{r-1} u \neq 0$.  This proves that $F u =
E_{j,i} u \neq 0$, so $E_{j,i}$ is injective, as desired.
\end{proof}

Consider the action of $\Sym_n$ (on the right) on $V^{\otimes
  n}$ by place permutation:
\begin{equation}
  (u_1 \otimes u_2 \otimes \cdots \otimes u_n)w = u_{1w^{-1}} \otimes
  u_{2w^{-1}} \otimes \cdots \otimes u_{nw^{-1}}
\end{equation}
for $w\in \Sym_n$, $u_1, \dots, u_n \in V$.  Note that $\Sym_n$ acts
(on the right) on the set $I(m,n)$, by the rule
\begin{equation}
  \mathbf{i}\cdot w = (i_{1w^{-1}}. i_{2w^{-1}}, \dots, i_{nw^{-1}})
\end{equation}
for $\mathbf{i}=(i_1, i_2, \dots, i_n)$, $w \in \Sym_n$. With this
notation, the action of $\Sym_n$ on the basis elements $e_\mathbf{i}$
of $V^{\otimes n}$ is given by $e_\mathbf{i} w = e_{\mathbf{i}\cdot
  w}$.

We now define a right action of $\He = \He_{\Q(v)}$ on $V^{\otimes n}$ on
basis elements by
\begin{equation}
e_{\mathbf{i}}T_{k} = \begin{cases}
v^2 e_{\mathbf{i}} & \text{ if $i_k = i_{k+1}$}\\
v e_{\mathbf{i}\cdot s_k} & \text { if $i_k < i_{k+1}$}\\
v e_{\mathbf{i}\cdot s_k} + (v^2 - 1)e_{\mathbf{i}}  & 
\text { if $i_k > i_{k+1}$.}
\end{cases}
\end{equation}
for any $\mathbf{i} \in I(m,n)$, and any $k = 1, \dots, n-1$. Notice
that when $v=1$ this is the same as the place permutation action of
$\Sym_n$. One may check that the actions of $\He$ and $\UU$ on
$V^{\otimes n}$ commute; this goes back to \cite{Jimbo}. The
commutativity of the two actions means that a weight space
$(V^{\otimes n})_\lambda$ is a right $\He$-module, for any $\lambda
\in X$.

\begin{lem}\label{lem:5.13}
  Let $\lambda \vDash n$, and assume $m \ge \ell(\lambda)$. Identify
  $\lambda = (\lambda_1, \dots, \lambda_m)$ with the element $\sum_i
  \lambda_i \varepsilon_i$ of $X$. Write $M^\lambda :=
  M_{\Q(v)}^\lambda$. Then $M^\lambda \simeq (V^{\otimes n})_\lambda$
  as right $\He$-modules.
\end{lem}

\begin{proof}
  (See \cite[\S 2]{DJ3} or \cite[\S 1]{GrLu}.)  Any $\underline{i} \in
  I(m,n)$ determines a unique row-standard tableau $t(\underline{i})$,
  in which $j$ appears in row $k$ whenever $i_j = k$, for $j=1,
  \ldots, n$.  For example,
  \[ 
  t((2,1,2,3))= \tiny\young(2,13,4)
  \]
  Let $\lambda \vDash n$ and let $\ell (\lambda) \leq m$.  Define a
  map $\phi : (V^{\otimes n})_\lambda \rightarrow M_R^\lambda$ on a
  basis by
  \[ 
   \phi(e_{\underline{i}}) = v^N x_{\lambda t(\underline{i})}, 
  \]
  where $\underline{i} \in I(m,n)$ indexes a simple tensor
  $e_{\underline{i}}$ of weight $\lambda$, and $N$ is the number of
  pairs $s < s'$ such that $i_s < i_{s'}$.  The map $\phi$ gives the
  desired isomorphism.
 \end{proof}

 The next result does not seem to have been previously observed in the
 literature. By specializing $v$ to 1 it gives in particular
 embeddings of permutation modules for $\Q\Sym_n$. 

\begin{lem}\label{lem-two}
  Suppose that $R = \Q(v)$ with $v$ an indeterminate, and $q=v^2$.  If
  $\lambda \dom \mu$, then $M^\lambda$ is isomorphic with an
  $\He$-submodule of $M^\mu$.
\end{lem}

\begin{proof}
  It suffices to prove the lemma for the case $\lambda \sdom \mu$ and
  $\lambda$, $\mu$ are adjacent in the dominance order.  In this case,
  the Young diagram of $\lambda$ is obtained from that of $\mu$ by
  raising one box from the $j^{th}$ to the $i^{th}$ row, where $i <
  j$.  This implies the statement in general, since whenever $\lambda
  \sdom \mu$, one gets from $\mu$ to $\lambda$ by a finite succession
  of such box-raising operations.
  
 Now, for $\lambda \sdom \mu$ and $\lambda$, $\mu$ adjacent, we have
 $\lambda = \mu + \varepsilon_i - \varepsilon_j$, and the result
 follows from Lemmas \ref{lem:Tony} and \ref{lem:5.13}.
\end{proof}

%%%%%%%%%%%%%%%%%%%%%%%%%%%%%%%%%%%%%%%%%%%%%%%%%%%%%%%%%%%%
\begin{rmk} \label{rmk:two} Suppose that $R$ is a field and $v=1$. It
  is easy to see that Lemma \ref{lem-two} may be false if the
  characteristic of $R$ is positive. For example, consider the map
  $f_{12}: M^{(3,1)} \to M^{(2,2)}$. The matrix of this map, with
  respect to a convenient choice of ordering of bases, is as follows:
  \[
  \left(\begin{smallmatrix}
  1&1&0&1&0&0\\
  1&0&1&0&1&0\\
  0&1&1&0&0&1\\
  0&0&0&1&1&1
  \end{smallmatrix}\right)
  \]
  and one observes that each column has precisely two 1's, so the sum
  of the rows is zero in characteristic 2. Thus, in characteristic 2
  the matrix has rank strictly less than 4, and thus the map is not
  injective. (This example is related to Example \ref{ex:counter}
  below.)

  It is easy to construct similar examples showing that Lemma
  \ref{lem-two} fails in any given positive characteristic, in the
  $v=1$ case.
\end{rmk}
%%%%%%%%%%%%%%%%%%%%%%%%%%%%%%%%%%%%%%%%%%%%%%%%%%%%%%%%%%%%%%%%%%

\section{The annihilator in the integral case}\label{sec:integral}
\noindent
In this section we work over the ring $\A=\Z[v,v^{-1}]$.  In Theorem
\ref{thm:main-integral} below, we describe $\ann_{\He_\A}
M_\A^\lambda$, using a refinement of an argument of \cite{Har}.  The
following result will be crucial in the proof of Theorem
\ref{thm:main-integral}.

\begin{lem}\label{lem-three}  If $\lambda \dom \mu$ for
$\lambda, \mu \vdash n$, then $\ann_{\He_\A} M_\A^\lambda \supseteq
  \ann_{\He_\A} M_\A^\mu$.
\end{lem}

\begin{proof}
This is an immediate consequence of Lemma \ref{lem-two}, since the
restriction of an injective map is injective.
\end{proof}

Next, we will need a result of Murphy.  Let $a,b \in \He_\A$, and let
$(a,b)$ denote the coefficient of $T_1$ in the expression $ab^{*} =
\sum_{w\in \Sym_n} c_w\, T_w$, where $c_w \in \A$. Then $(\ ,\ )$ is a
non-degenerate, symmetric bilinear form on $\He_\A$. This bilinear
form satisfies the properties
\begin{equation} \label{eqn:bilform}
  (a, bd) = (ad^{*}, b); \quad (a, db) = (d^{*}a, b)
\end{equation}
for any $a,b,d \in \He_\A$.

\begin{lem}[Murphy {\cite[Lemma 4.16]{Murphy2}}]\label{lem:keylem} 
Let $s$ and $t$ be row-standard $\mu$-tableaux and let $u$,$w \in
\Tab(\lambda)$, where $\mu \vDash n$ and $\lambda \vdash n$. Then:

(a) $(x_{st}, x^\sharp_{uw}) = 0$ unless $(\tr{u},\tr{w}) \dom (s,t)$;

(b) $(x_{\tr{u}\tr{w}}, x^\sharp_{uw}) = \pm v^{2b}$ where $b=
\ell(d(t_\lambda)) + \sum_{i \geq 1}\lambda'_i(\lambda'_i -1)/2$.
\end{lem}

The following result is an ``integral'' version of Theorem
\ref{thm:main-ss}.

\begin{thm}\label{thm:main-integral}
   For any $\lambda \vdash n$ we have $\ann_{\He_\A} M_\A^\lambda =
   {\He_\A}[\notlessdom \tr{\lambda}]^\sharp$.
\end{thm}

\begin{proof}
By Lemma \ref{lem-one}, we have already the containment
${\He_\A}[\notlessdom \tr{\lambda}]^\sharp \subseteq \ann_{\He_\A}
M_\A^\lambda$, so we have only to prove the reverse containment.  Let
\[
0 \ne a = \mathop{\sum_{s,t \in \Tab(\lambda)}}_{\lambda \vdash n}
a_{st}x^\sharp_{st} \in \ann_{\He_\A} M_\A^\lambda,
\] 
where $a_{st} \in \A$.  It suffices to prove that $a_{st} = 0$ for all
$s$ and $t$ of some shape $\mu \lessdom \tr{\lambda}$.

Suppose not. By Lemma \ref{lem-one} we have $\sum_{(s,t) \in \Phi}
a_{st}x^\sharp_{st} \in \ann_{\He_\A} M_\A^\lambda$, hence it follows
that
\[
  0 \ne a_0 = \sum_{(s,t) \in \comp{\Phi}} a_{st}x^\sharp_{st} \in
  \ann_{\He_\A} M_\A^\lambda,
\]
where $\Phi = \bigsqcup_{\mu \notlessdom \tr{\lambda}} \Tab(\mu)
\times \Tab(\mu)$ and $\comp{\Phi} = \bigsqcup_{\mu \lessdom
  \tr{\lambda}} \Tab(\mu) \times \Tab(\mu)$.  Let $(s_0,t_0)$ be a
minimal pair in $\comp{\Phi}$ such that $a_{s_0t_0} \ne 0$; i.e.,
\[
a_{st} = 0 \text{ for all $(s,t) \in \comp{\Phi}$ satisfying
$(s, t) \slessdom (s_0,t_0)$} .
\]
Let $\lambda_0$ be the shape of $\tr{s_0}$ (= shape of
$\tr{t_0}$). Then $\tr{\lambda_0}$ is the shape of $s_0, t_0$, so
$\tr{\lambda_0} \lessdom \tr{\lambda}$, and hence $\lambda \lessdom
\lambda_0$. By Lemma \ref{lem-three}, anything annihilating
$M_\A^\lambda$ also annihilates $M_\A^{\lambda_0}$. Thus, it follows that
$x_{\lambda_0} T_{d(\tr{t_0})} a_0 = 0$. Hence,
$T^*_{d(\tr{s_0})}x_{\lambda_0} T_{d(\tr{t_0})} a_0 =
x_{\tr{s_0}\tr{t_0}}\, a_0 = 0$. So, by the definition of the bilinear
form, we have
\begin{equation*}
  0 = \left( x_{\tr{s_0}\tr{t_0}} , \sum_{(s,t) \in \comp{\Phi}}
  a_{st}x^\sharp_{st} \right) = \sum_{(s,t) \in \comp{\Phi}} a_{st}
  (x_{\tr{s_0}\tr{t_0}}, x^\sharp_{st}).
\end{equation*}
By Lemma \ref{lem:keylem}(a), all the terms in the last sum are zero
unless $(\tr{s}, \tr{t}) \dom (\tr{s_0}, \tr{t_0})$, i.e.,\ unless
$(s,t) \lessdom (s_0,t_0)$.  By the minimality assumption, $a_{st}=0$
for all pairs $(s,t)$ strictly less dominant than $(s_0,t_0)$. Thus,
the sum $\sum_{(s,t) \in \comp{\Phi}} a_{st} (x_{\tr{s_0}\tr{t_0}},
x^\sharp_{st})$ collapses to a single term $a_{s_0t_0} (
x_{\tr{s_0}\tr{t_0}}, x^\sharp_{s_0 t_0} )$, and by our assumption and
Lemma \ref{lem:keylem}(b) this term is nonzero. This is a
contradiction.  This contradiction establishes the desired opposite
inclusion, and proves the theorem.
\end{proof}

\begin{ex}\label{ex:counter}
We can now give an example to show that the annihilator of
$M_R^\lambda$ depends on $R$, even when $R$ is a field. We take $v=1$,
$\lambda=(2,2)$, and let $R$ be a field of characteristic 2.  It is
quickly seen that the element $r=(23)+(1342)+(1243)+(14)$
% r=(3124)+(3142) + (2413)+(4231)$ 
(written in cycle notation) annihilates $M_R^{(2,2)}$ in
characteristic 2.
%since
%\begin{eqnarray*}
%(1,1,2,2) r & = & (1,2,1,2) + (1,2,1,2) + (2,1,2,1) + (2,1,2,1) = 0\\
%(1,2,1,2)r & = & (1,1,2,2) +(2,2,1,1) +(1,1,2,2) + (2,2,1,1) = 0\\
%(1,2,2,1)r &=& (1,2,2,1) + (2,1,1,2) + (2,1,1,2) + (1,2,2,1) =0\\
%(2,1,1,2)r&=& (2,1,1,2) + (1,2,2,1) + (1,2,2,1) + (2,1,1,2) = 0\\
%(2,1,2,1)r & = & (2,2,1,1) + (1,1,2,2) + (2,2,1,1) + (1,1,2,2) =0 \\
%(2,2,1,1)r& = & (2,1,2,1) + (2,1,2,1) +(1,2,1,2) + (1,2,1,2) = 0.
%\label{ex:char2ann}
%\end{eqnarray*}
However, we claim $r$ is not in the span of the basis elements of
$\ann_{\Z\Sym_4}M_\Z^{(2,2)}$ with coefficients reduced modulo 2.
 
By Theorem \ref{thm:main-integral}, for standard tableaux $a=
\tiny\young(1234)$, $b= \tiny\young(123,4)$, $c=\tiny\young(124,3)$,
and $d=\tiny\young(134,2)$, the following 10 Murphy elements form a basis
for $\ann_{\Z\Sym_4}M_\Z^{(2,2)}$.
\begin{eqnarray*}
x^\sharp_{aa}&=& \sum_{w \in \Sym_4}(-1)^{\ell(w)} w \\
x^\sharp_{bb} & = & (1) - (12)-(13)-(23) + (123)+(132)\\
x^\sharp_{bc} & = & (34)-(12)(34)-(143)-(243)+(1243)+(1432)\\
x^\sharp_{bd}&=&(234)-(1342)-(1423)-(24)+(13)(24)+(142)\\
x^\sharp_{cc}&=& (1)-(12)-(14)-(24)+(124) +(142)\\
x^\sharp_{cb} & = &(34)-(12)(34)-(134)-(234)+(1234)+(1342) \\
x^\sharp_{cd}& = &(23)-(132)-(14)(23)-(243)+(1324)+(1432) \\
x^\sharp_{db}& = & (243)-(1243)-(1324) - (24) + (124) + (13)(24)\\
x^\sharp_{dc} & = & (23) - (123) - (14)(23) - (234) + (1234) + (1423) \\
x^\sharp_{dd} & = & (1) - (13) - (14) - (34) +(134) + (143).
\label{ex:basis22}
\end{eqnarray*}
It is easy to see that $r$ does not belong to the $R$-linear span of
these basis elements.
%
%Suppose $r$ can be written as the sum of basis elements with
%coefficients reduced modulo 2.  We consider the cases in which
%$x^\sharp_{aa}$ is and is not present in this sum. If $x^\sharp_{aa}$
%does not occur, one of $x^\sharp_{cc}$ or $x^\sharp_{dd}$ must be
%present, since (14) appears in $r$.  If $x^\sharp_{cc}$ appears,
%$x^\sharp_{bb}$ must also appear, since (12) occurs only in those two
%(non-$x^\sharp_{aa}$) basis elements and does not appear in $r$.
%Since (13) appears only in $x^\sharp_{bb}$ and $x^\sharp_{dd}$, the
%presence of $x^\sharp_{bb}$ implies the presence of $x^\sharp_{dd}$.
%This in turn results in (14) occurring an even number of times, in
%which case it disappears from $r$.  The case in which (14) is
%represented by $x^\sharp_{dd}$ is similar.
%
%Next, suppose $x^\sharp_{aa}$ appears in the sum.  Since (14)(23) is
%not in $r$ and, apart from $x^\sharp_{aa}$, appears only in
%$x^\sharp_{cd}$ and $x^\sharp_{dc}$, one, but not both, of
%$x^\sharp_{cd}$ or $x^\sharp_{dc}$ must be present.  For (23) to
%appear in $r$, two of $x^\sharp_{bb}$, $x^\sharp_{cd}$ or
%$x^\sharp_{dc}$ must be in the sum.  However, $x^\sharp_{bb}$ and
%$x^\sharp_{cd}$ both appearing would imply (132) occurs an odd number
%of times, and $x^\sharp_{bb}$ and $x^\sharp_{dc}$ both appearing would
%imply (123) occurs an odd number of times.  Since neither (123) nor
%(132) appear in $r$, there is no combination of basis elements with
%coefficients reduced modulo 2 which gives $r$.
%
This shows that $\dim_R \ann_{R\Sym_4} M_R^{(2,2)} \ge 11$. (In fact,
it equals 11.)
\end{ex}

\section{Applications}\noindent
We give some consequences of Theorem \ref{thm:main-integral}. As
above, we work over $\A= \Z[v, v^{-1}]$ in this section. First we note
the following consequence of the main results.

\begin{rmk}
  As already noted in Section 3, we obtain the following result
  immediately by applying the involution $\sharp$ to the equalities in
  Theorems \ref{thm:main-ss} and \ref{thm:main-integral}:
  \[ 
   \ann_{\He_R} \tilM_R^\lambda = {\He_R}[\notlessdom \tr{\lambda}]  
  \]
  provided ${\He_R}$ is semisimple over a field $R$, or $R = \A$.
\end{rmk}

The next application is the cellularity of the algebras
${\He_\A}/(\ann_{\He_\A} M_\A^\lambda)$.

\begin{cor}
  Let $\lambda \vdash n$.  The quotient algebra
  ${\He_\A}/(\ann_{\He_\A} M_\A^\lambda)$ is cellular with cell basis
  $\{ x^\sharp_{st} + {\He_\A}[\notlessdom \tr{\lambda}]^\sharp : s,t
  \in \Tab(\mu), \mu \lessdom \tr{\lambda} \}$.
\end{cor}

\begin{proof}
  This follows immediately from the Theorem and the theory of cellular
  algebras.
\end{proof}

Finally, we observe that certain other modules have the same
annihilator as $M_\A^\lambda$. Recall that ${\He_\A}[\dom
\lambda]^\sharp/{\He_\A}[\sdom \lambda]^\sharp$ is naturally a right
${\He_\A}$-module under right multiplication by elements of ${\He_\A}$
(see \S\ref{sec:Murphy}).

\begin{cor} For any $\lambda \vdash n$, we have the equalities
\begin{align*}  
\ann_{\He_\A} {\He_\A}[\dom \lambda] &= \ann_{\He_\A}
  ({\He_\A}[\dom \lambda]/{\He_\A}[\sdom \lambda])\\ 
  &= \ann_{\He_\A} M_\A^\lambda 
   = {\He_\A}[\notlessdom \tr{\lambda}]^\sharp .
\end{align*} 
\end{cor}

\begin{proof}
  Recall that $M_\A^\lambda$ has a basis given by all $x_\lambda
  T_{d(t)}$ as $t$ ranges over the set of row-standard tableaux of
  shape $\lambda$.  Now, any $a \in {\He_\A}$ acting as zero on
  $x_\lambda T_{d(t)}$ also acts as zero on $T^*_{d(s)} x_\lambda
  T_{d(t)}$, for any $s$. Thus, if $a \in \ann_{\He_\A} M_\A^\lambda$,
  then $a$ acts as zero on all $x_{st} = T^*_{d(s)} x_\lambda
  T_{d(t)}$, with $s,t \in \Tab(\lambda)$.  Since the images of these
  elements in ${\He_\A}[\dom \lambda]/{\He_\A}[\sdom \lambda]$ form a
  basis for that ${\He_\A}$-module, it follows that $a \in
  \ann_{\He_\A} (A^x[\dom \lambda]/A^x[\sdom \lambda])$. This proves
  the inclusion $\ann_{\He_\A} M_\A^\lambda \subseteq \ann_{\He_\A}
  ({\He_\A}[\dom \lambda]/{\He_\A}[\sdom \lambda])$. The opposite
  inclusion is clear, as $M_\A^\lambda$ is isomorphic with a submodule of
  ${\He_\A}[\dom \lambda]/{\He_\A}[\sdom \lambda]$. This proves the
  second equality.

  By Lemma \ref{lem-three}, we know that $\ann_{\He_\A} M_\A^\lambda
  \subseteq \ann_{\He_\A} M_\A^\mu$ for any $\lambda \lessdom \mu$. Thus,
  by the result of the previous paragraph, it follows that any $a \in
  \ann_{\He_\A} M_\A^\lambda$ acts as zero on all $x_{uw}$ for $(u,w)$ of
  shape $\mu$, for any $\mu \dom \lambda$. Since ${\He_\A}[\dom
    \lambda]$ is generated over $\A$ by such $x_{uw}$, it follows that
  $\ann_{\He_\A} M_\A^\lambda \subseteq \ann_{\He_\A} {\He_\A}[\dom
    \lambda]$. The opposite inclusion is clear, so this proves the
  first equality. The third equality is from Theorem
  \ref{thm:main-integral}.
\end{proof}

\bibliographystyle{amsalpha}

\begin{thebibliography}{100}\frenchspacing\raggedright
\small %Use smaller font size

\bibitem{DJ:Hecke} R.~Dipper and G.~James, Representations of
  Hecke algebras of general linear groups, \emph{Proc. London
    Math. Soc.} (3) \textbf{52} (1986), 20--52.

\bibitem{DJ:qsa} R.~Dipper and G.~James, The $q$-Schur
  algebra, \emph{Proc. London Math. Soc.} (3) \textbf{59} (1989),
  23--50.

\bibitem{DJ3} R.~Dipper and G.~James, $q$-tensor space and
  $q$-Weyl modules, \emph{Trans. Amer. Math. Soc.} (1) \textbf{327}
  (1991), 251--282.

\bibitem{DJ4} R.~Dipper and G.~James, Blocks and idempotents
  of Hecke algebras of general linear groups, \emph{Proc. London
    Math. Soc.} (3) \textbf{54} (1987), 57--82.

\bibitem{Tony} A. Giaquinto, Quantization of tensor
  representations and deformation of matrix bialgebras, {\em J.
  Pure Appl. Algebra}, \textbf{79} (1992), 169--190.


\bibitem{GL} J.J. Graham and G.I. Lehrer, Cellular algebras,
{\em Invent. Math.} {\bf123} (1996), 1--34.

\bibitem{GrLu} I.~Grojnowski and G.~Lusztig, On bases of
  irreducible representations of quantum $GL_n$, in ``Kazhdan-Lusztig
  theory and related topics'', {\em Contemp. Math.} {\bf 139}
  (1992), 167--174.

\bibitem{Har} M. H\"{a}rterich, Murphy bases of
 generalized Temperley-Lieb algebras, {\em Archiv Math.} {\bf 72}
 (1999), 337--345.

\bibitem{James} G.D. James, {\em The Representation Theory of
the Symmetric Groups}, Lecture Notes in Math. {\bf682},
Springer-Verlag, Berlin 1978.

\bibitem{Jimbo} M. Jimbo, A $q$-analogue of $U(\mathfrak{gl}
  (N+1))$, Hecke algebra, and the Yang-Baxter equation.
  \textit{Lett. Math. Phys.} \textbf{11} (1986), 247--252.

\bibitem{Mathas} A.~Mathas, \emph{Iwahori--Hecke algebras and
  Schur algebras of the symmetric group}, University Lecture Series,
  15, American Mathematical Society, Providence, RI, 1999.

\bibitem{Murphy} G.E. Murphy, On the representation theory of
the symmetric groups and associated Hecke algebras, {\em J. Algebra}
{\bf 152} (1992), 492--513.

\bibitem{Murphy2} G.E. Murphy, The representations of Hecke
  algebras of type $A_n$, {\em J. Algebra} {\bf 173} (1995), 97--121.

\bibitem{Xi} N.~Xi, Root vectors in quantum groups,
  \emph{Comment. Math. Helv.}  \textbf{69} (1994), 612--639.

\end{thebibliography}

\end{document}